\documentclass{article}
\usepackage{amsthm,amsmath,amssymb}
\usepackage[colorlinks,citecolor=blue,urlcolor=blue]{hyperref}  
\usepackage{enumerate}

\newtheorem{theorem}{Theorem}[section]

\newtheorem{corollary}[theorem]{Corollary} 
\theoremstyle{definition}
\newtheorem{definition}[theorem]{Definition}

\theoremstyle{remark}
\newtheorem{remark}[theorem]{Remark}
\newtheorem*{note}{Note}
\newtheorem{example}[theorem]{Example}
\newtheorem{proposition}{Proposition}[theorem]

\date{August 19, 2014}

\newcommand{\J}{J}
\newcommand{\N}{\mathbb{N}}
\newcommand{\p}{\mathsf{p}}
\renewcommand{\k}{\mathsf{k}}
\newcommand{\s}{\mathsf{s}}
\newcommand{\true}{\mathsf{T}}
\newcommand{\false}{\mathsf{F}}
\newcommand{\halts}{\downarrow}
\newcommand{\ass}[1]{\text{Ass}(#1)}
\newcommand{\rt}[1]{\text{RT}(#1)}

\newcommand{\la}{\langle}
\newcommand{\ra}{\rangle}

\newcommand{\varco}{\thinspace | \thinspace}

\title{Effective Operations of Type 2 in Pcas}
\author{Eric Faber \and
Jaap van Oosten \\
Department of Mathematics,
    Utrecht University\\
    P.O. Box 80.010, 3508TA Utrecht \\
    The Netherlands\\
J.vanOosten@uu.nl}
\begin{document}
\maketitle
\begin{abstract}
We exhibit a way of ``forcing a functional to be an effective operation'' for arbitrary partial combinatory algebras (pcas). This gives a method of defining new pcas from old ones for some fixed functional, where the new partial functions can be viewed as computable relative to that functional. It is shown that this generalizes a notion of computation relative to a functional as defined by Kleene for the natural numbers.

The construction can be used to study subtoposes of the Effective Topos.
We will do this for a particular functional that forces every arithmetical set to be decidable.

In this paper we also prove the convenient result that the two definitions of a pca that are common in the literature are essentially the same.
\end{abstract}

\section{Introduction}
Effective operations have originally been studied in computability theory, with as famous results the theorems of Myhill-Shepherdson and Rice-Shapiro.
They have also been studied in relation to computable functionals (e.g. see \cite{HYLAND-GANDY}). For example, it is known that
not every functional that is an effective operation is a computable functional.
In realizability, partial combinatory algebras (pcas) are used as models for abstract Turing machines.
The notion of effective operation easily extends to pcas.
In this paper, we discuss effective operations on arbitrary partial combinatory algebras, and present a way to force a functional
to be an effective operation: for a (total) functional $F: A^A \rightarrow A$, where $A$ is a pca, we construct a pca $A[F]$ for which $F$ is an effective operation.

In \cite{KLEENE-1959}, Kleene introduced the notion of a \emph{computable functional}, and with it the notion of computability relative to a functional.
It turns out that for $A = \mathcal{K}_1$, Kleene's first model, the functions computable in $\mathcal{K}_1[F]$ are precisely the functions computable
relative to $F$ in Kleene's sense. This yields an alternative approach to computation in a type 2 functional that generalizes to arbitrary pcas.
We can use this to study local operators in realizability toposes.
The theory of recursive hierarchies of functionals (e.g. see \cite{HINMAN}) can therefore be applied to study subtoposes of realizability toposes.

As an example, we look at a specific local operator $\J$ that is studied  in \cite{OOSTEN:arXiv1301.0735}.
It turns out that the functions ``computable in $\J$'' are precisely the functions computable relative to a certain functional.
On the way, we also prove that a strengthening of our definition of a pca, that we call strict pcas, is (up to isomorphism) not a proper strengthening at all.
This is posed as an open question in \cite{FREY}.
The result is convenient, since in the literature there has not been wide consensus on the definition of pca; several authors have been using strict pcas by default.

\section{Preliminaries}
We start by reviewing the definition of a partial combinatory algebra and some basic constructions.
\begin{definition}
 A partial applicative structure on a set $A$ is a partial map $A \times A \rightarrow A$ that we denote by
 \[
  (a,b) \mapsto ab.
 \]
\end{definition}
We write $ab \halts$ to say that $ab$ is defined, $ab = c$ to say that $ab$ is defined and has value $c$.
In writing compositions of the application, we adopt the convention of \emph{associating to the left}, so:
\[
 abcd = (((ab)c)d).
\]

For $V = \{x_0,x_1,\ldots\}$ an enumerable set of variables, we define the set of \emph{terms} by:
\begin{definition}
For all $a \in A$, $x \in V$
\begin{enumerate}
 \item $a$ is a term
 \item $x$ is a term
 \item For $t,s$ terms, $(ts)$ is a term.
\end{enumerate}
\end{definition}
A term $t$ without free variables is called \emph{closed} and might \emph{denote} or not. This is defined inductively:
\begin{enumerate}
 \item $a$ denotes and has value $a$
 \item If $t$ denotes and has value $a$, and $s$ denotes and has value $b$, and $ab = c$,
 then $(ts)$ denotes and has value $c$.
\end{enumerate}
We write $t \halts$ if $t$ denotes, and $t = a$ if $t$ denotes and has value $a$.
We write $t \lesssim s$ if $s \halts$ implies that $t \halts$, and in that case $s = t$.
We write $t \simeq s$ if $t \lesssim s$ and $s \lesssim t$.

We denote a term $t$ with free variables $x_1,\ldots,x_n$ by \[t(x_1,\ldots,x_n)\].
For $a_1,\ldots, a_n \in A$, we can substitute $x_1,\ldots,x_n$ by $a_1,\ldots,a_n$,
resulting in a closed term $t(a_1,\ldots,a_n)$.

\begin{definition}\label{dfn:pca}
 A partial combinatory algebra is a partial applicative structure $A$, with elements $\k,\s \in A$ such that for all $a,b \in A$:
 \begin{enumerate}
  \item $\s a b \halts$ 
  \item $\k a b = a$ 
  \item $\s a b c \lesssim a c (b c)$.
 \end{enumerate}
 
 We also say that the applicative structure on $A$ is \emph{combinatory complete}.
\end{definition}

\begin{example}
 There are several examples of pcas, \cite{OOSTEN} treats a lot of them.
 Here we just list a few with little explanation:
 \begin{itemize}
  \item The trivial pca, $A = \{*\}$ with $** = *$.
  \item \emph{Kleene's first model} $\mathcal{K}_1$, with underlying set $\N$.
  Application is given by:
  \[
   nm \simeq \phi_n(m)
  \]
where $\phi_n$ is the $n-$th partial function for some model of computation and coding.
 \end{itemize}
\end{example}

An introduction to partial combinatory algebras and its basic properties can be found in chapter 1 of \cite{OOSTEN}.
Here we will just quickly present some tools that we need later on, the proofs can all be found in \cite{OOSTEN}.
For a pca $A$, an element $a \in A$ defines a partial function $A \rightarrow A$ by:
\[
 b \mapsto ab.
\]
We say that $a$ is an \emph{index} for this partial function.
By abuse of language, we sometimes say that $a$ is a partial (computable) function.

An important result for pcas is the \emph{recursion theorem}:
\begin{theorem}[Recursion Theorem]
 For every $f \in A$, there exists $e$ such that for all $a \in A$:
 \[
  ea \lesssim fea.
 \]
\end{theorem}

Also important is the fact that a pca is \emph{combinatory complete}.
This means that for every term $t(x_1,\ldots,x_n)$ with free variables $x_1,\ldots,x_n$,
there is an element $\la x_1\cdots x_n \ra t \in A$ such that for all $a_1,\ldots,a_n \in A$:
\begin{align}
 (\la x_1\cdots x_n \ra t)a_1\cdots a_{n-1} \halts \\
  (\la x_1\cdots x_n \ra t)a_1\cdots a_{n-1} a_n \lesssim t(a_1,\ldots,a_n).
\end{align}
The notation $\la x_1\cdots x_n \ra t$ should remind the reader of $\lambda-$abstraction in $\lambda-$calculus.
We use a different notation because it is slightly different from $\lambda-$abstraction, but it is applied in the same way.

Besides $\k$ and $\s$, we identify (a choice of) other canonical elements of $A$, among which:
\begin{itemize}
 \item Boolean combinators $\true$ and $\false$ that satisfy for all $a,b \in A$:
 \begin{align}
  \true a b = a \\
  \false a b = b.
 \end{align}
Indeed, one can take $\true = \k$.
\item Using combinatory completeness and the boolean combinators, we can define if/else statements:
for $t,s$ closed terms:
\[
 \left(\la v \ra v (\la x \ra t)(\la x \ra s)\k\right) b \lesssim \begin{cases}
                                                  t &\text{ if } b = \true \\
                                                  s &\text{ if } b = \false \\
                                                  \text{unspecified} &\text{ otherwise}.
                                                 \end{cases}
\]
For legibility, we denote a term like $v (\la x \ra t)(\la x \ra s)\k$ by
\[
 \textsf{if } v \textsf{ then } t \textsf{ else } s.
\]

\item The pairing combinator $\p$ together with projections $\p_0,\p_1$:
\begin{align}
 \p_0(\p ab) = a \\
 \p_1(\p ab) = b
\end{align}
\item The Curry numerals: each pca contains representables of the natural numbers, i.e. for any $n \in \N$ an element $\overline{n} \in A$. It is a basic result (for a proof, see \cite{OOSTEN}) that, for every $k$, there is for every $k-$ary partial recursive function $F: \N^{k} \rightarrow \N$ an index $f \in A$ such that for all $(n_1,\ldots,n_k) \in \operatorname{dom}{F}$:
\[
	f\overline{n_1}\dots\overline{n_k} = \overline{F(n_1,\ldots,n_k)}.
\]
\end{itemize}

When there is chance of confusion, we decorate (choices of) canonical elements by a subscript to denote the pca to which it belongs, e.g. $\k_A,\s_A,\true_A$,etc.

We can also code tuples of elements in a computable way.
That is, there is an index $t \in A$ such that for all $n \in \N$, $u_0,\ldots,u_n \in A$ we can code a tuple
\[
 [u_0,\ldots,u_n] := t\overline{n+1}u_0\cdots u_n \in A
\]
such that projection, concatenation of tuples, and a length function (that yields the length of a tuple as a curry numeral) are all computable.
We denote concatenation by $*$:
\[
 [u_0\ldots,u_i]*[u_{i+1},\ldots,u_n] = [u_0,\ldots,u_n]
\]
and for $u = [u_0,\ldots,u_n]$, we write $u^{<{i+1}}$ for $[u_0,\ldots,u_i]$.
These basic operations can all be done computably and uniformly in the variables.

There is a notion of morphism between partial combinatory algebras due to John Longley (see \cite{LONGLEY}), called \emph{applicative morphism}.
In the following definition, we adopt the following notation:
for a total relation $R: X \rightarrow Y$ between sets $X$ and $Y$, we denote for every $x \in X$ by $R(x)$ the set:
\[
 \{y \in Y \varco xRy\}.
\]
For $A$ a pca, $a \in A$ and $P \subseteq A$ a subset, we write $aP \halts$ if for all $b \in P$, $ab\halts$,
and in that case
\[
 aP := \{ab \varco b \in P\}.
\]
For $P,Q \subseteq A$, we write $PQ \halts$ if for all $b \in P, b' \in Q$, $bb' \halts$, and in that case
\[
 PQ = \{bb' \varco b \in P, b' \in Q\}.
\]
Whenever we write compositions of these applications of an element and a subset, or a subset and a subset, we again associate to the left.
\begin{definition}[Longley]
 Let $A,B$ be pcas. An applicative morphism $\gamma: A \rightarrow B$ is a total relation from $A$ to $B$, together with an element $r \in B$
 such that whenever $a, a' \in A$ and $aa'\halts$, we have:
 \[
   r\gamma(a)\gamma(a') \halts \text{ and } r\gamma(a)\gamma(a') \subseteq \gamma(aa'). 
 \]
The element $r$ is said to \emph{realize} $\gamma$.

An applicative morphism $\gamma$ is \emph{decidable} if there is $d \in B$ such that
\begin{align}
 d\gamma(\true_A) = \{\true_B\} \\
 d\gamma(\false_A) = \{\false_B\}.
\end{align}
We call $d$ a \emph{decider} for $\gamma$.
\end{definition}
The composition of two applicative morphisms is defined by composition of the relations, and one can verify that this is again an applicative morphism.
There is a pre-order on applicative morphisms defined as follows: for $\gamma: A \rightarrow B$, $\delta: A \rightarrow B$,
write $\gamma \preceq \delta$ when there is $t \in B$ such that for all $a \in A$:
\[
 t\gamma(a) \subseteq \delta(a).
\]
We write $\gamma \simeq \delta$ if $\gamma \preceq \delta$ and $\delta \preceq \gamma$, in this case $\gamma$ and $\delta$ are \emph{equivalent}.
One can show that this pre-order is preserved by composition on both sides.
We thus obtain a pre-order enriched category PCA consisting of partial combinatory algebras as objects, and applicative morphisms as arrows.
The identity arrow for a pca $A$ is denoted by $\iota_A$ and given by the identity relation on $A$.
\begin{note}
	As remarked in the introduction, definition \ref{dfn:pca} is weaker than the definition that is used in some other literature (including \cite{OOSTEN}).
	However, it will be a consequence of theorem \ref{thm:pca_iso_strict_pca} below that the corresponding categories of pcas, with morphisms as above, are equivalent.
\end{note}

For a pca $A$, we define a set $\text{PR}_A$ by:
\[
 \text{PR}_A = \{f: A \rightarrow A \text{ partial }\varco (\exists a \in A)(\forall b \in A) ab \lesssim f(b) \}.
\]
We also define a function $I_1^A: \text{PR}_A \rightarrow A$ by:
\[
 I_1^A(f) = \{a \in A \varco (\forall b \in A) ab \lesssim f(b)\}.
\]
For $\gamma: A \rightarrow B$ an applicative morphism, we define a set $\text{PR}_\gamma$ by:
\begin{align}
 \text{PR}_\gamma: \{f: A \rightarrow A \text{ partial } \varco & (\exists b \in B)(\forall a \in A) f(a) \halts \Rightarrow b\gamma(a) \subseteq \gamma(f(a)) \}.
\end{align}
Lastly, we also define a function $I_1^\gamma: \text{PR}_\gamma \rightarrow B$ by:
\[
 I_1^\gamma(f) = \{b \in B \varco (\forall a \in A) f(a) \halts \Rightarrow b\gamma(a) \subseteq \gamma(f(a)) \}.
\]
Observe that $\text{PR}_A = \text{PR}_{\iota_A}$.

The following result by the second author (see \cite{OOSTEN-GENERAL_RECURSION} or \cite{OOSTEN}) tells us that we can adjoin any partial function $f: A \rightarrow A$ to a pca
 $A$, to obtain a pca $A[f]$ in which $f$ has an index. This generalizes the notion of a Turing machine with oracle from ordinary computability theory to pcas.
 \begin{theorem}[\cite{OOSTEN-GENERAL_RECURSION}]\label{thm:general_recursion}
  For $A$ a pca and $f: A \rightarrow A$ a partial function, there exists a pca $A[f]$ defined on the set $A$ such that
  the identity on $A$ is a decidable applicative morphism $\iota_f: A \rightarrow A[f]$, and:
  \begin{enumerate}
   \item $f \in \text{PR}_{A[f]} = \text{PR}_{\iota_f}$.
   \item Whenever $\gamma: A \rightarrow B$ is a decidable applicative morphism such that $f \in \text{PR}_\gamma$,
   there exists a decidable $\gamma_f: A[f] \rightarrow B$ such that $\gamma_f \circ \iota_f = \gamma$.
   Moreover, if $\delta: A[f] \rightarrow B$ is a decidable applicative morphism such that $\delta \circ \iota_f \cong \gamma$, then $\delta \cong \gamma_f$.
  \end{enumerate}
 \end{theorem}
 
 The proof of the above theorem is done by defining an applicative structure on $A$ using \emph{dialogues}:
 \begin{definition}
  Let $f: A \rightarrow A$ be a partial function, and $a,b \in A$.
  An $f-$\emph{dialogue} between $a$ and $b$ is an element $u = [u_0,\ldots,u_n] \in A$ such that
  for all $i \leq n$:
  \[
   a([b]*u^{<i}) = \p \false v_i
  \]
where $f(v_i)$ is defined and $f(v_i) = u_i$.

An $f-$dialogue $u$ between $a$ and $b$ \emph{halts} if there exists $c \in A$ such that
\[
 a([b]*u) = \p \true c.
\]
We write $a \cdot^f b \halts$ if there is a (necessarily unique) halting $f-$dialogue $u$ between $a$ and $b$, and $a \cdot^f b = c$ if $a \cdot^f b \halts$
and
\[
 a([b]*u) = \p \true c
\]
where $u$ is a a halting $f-$dialogue between $a$ and $b$.
 \end{definition}

 The applicative structure on $A[f]$ is then given by $(a,b) \mapsto a\cdot^f b$.
 For a full proof of combinatory completeness, see \cite{OOSTEN} or \cite{OOSTEN-GENERAL_RECURSION}. Note that these sources use the strict version of a pca (definition \ref{dfn:strict_pca}), but the proofs work without any modification (alternatively, one can apply theorem \ref{thm:pca_iso_strict_pca}).
 The following corollary corresponds to corollary 1.3 in \cite{OOSTEN-GENERAL_RECURSION}:
 \begin{corollary}\label{cor:oosten_original}
 Let $A$ be a pca.
 
  \begin{enumerate}
   \item If $f \in \text{PR}_A$, then $A[f]$ and $A$ are isomorphic pcas.
   \item If $f, g: A \rightarrow A$ are partial functions, then $A[f][g]$ and $A[g][f]$ are isomorphic.
   \item If $\mathcal{K}_1^f$ is the pca of partial recursive application with an oracle for $f$,
   then $\mathcal{K}_1^f$ is isomorphic to $\mathcal{K}_1[f]$.
  \end{enumerate}
 \end{corollary}

 \section{Effective operations}
 We will now define the set of effective operations on pcas.
 For any pca $A$, let $\text{Tot}_A$ be the subset of $\text{PR}_A$ consisting of all total functions:
 \[
  \text{Tot}_A = \{ f \in \text{PR}_A \varco f \text{ is total }\}.
 \]
Similarly, define for any $\gamma$ the subset $\text{Tot}_\gamma \subseteq \text{PR}_\gamma$ consisting of all total functions in $\text{PR}_\gamma$.
\begin{definition}[Effective Operation]
For a pca $A$, we define the set $E_A$ of \emph{effective operations} in $A$ by:
\begin{align}
 E_A = \{F: A^A \rightarrow A \text{ partial} \varco (\exists a \in A)(\forall f \in \text{Tot}_A)F(f) \text{ defined } \Rightarrow \\
 (a I_1(f)\halts \text{ and } a I_1^A(f) = \{F(f)\})\}.
\end{align}
Define a function $I_2^A: E_A \rightarrow A$ by:
\[
 I_2^A(F) = \{a \varco (\forall f \in \text{Tot}_A)(\forall b \in I_1^A(f))(ab \lesssim F(f))\}
\]

Similarly, for $\gamma: A \rightarrow B$ an applicative morphism we define the set $E_\gamma$ of \emph{effective operations relative to} $\gamma$ by:
\begin{align}
 E_\gamma = \{F: A^A \rightarrow A \text{ partial} \varco (\exists b \in B)(\forall f \in \text{Tot}_\gamma)
 F(f) \text{ defined } \Rightarrow \\
 (b I_1^\gamma(f) \halts \text{ and } b I_1^\gamma(f) \subseteq \gamma(F(f)))\}.
\end{align}
\end{definition}

\begin{example}\label{ex:functionalE} As an example, we look at the following functional $E: A^A \rightarrow A$, where $A$ is an arbitrary pca.
 \[
  E(\alpha) = \begin{cases}
               \true &\text{ if } (\exists a) \alpha(a) = \true \\
               \false &\text{ otherwise}.
              \end{cases}
 \]
For $A = \mathcal{K}_1$, it is easily seen that the above is not an effective operation (otherwise the halting problem is decidable).
However, the below construction shows that there are a lot of pcas in which $E$ is an effective operation.
\end{example}

The statement of the following theorem is very similar to theorem \ref{thm:general_recursion}.
\begin{theorem}\label{thm:force_eff_operation}
 Let $A$ be a pca, and $F: A^A \rightarrow A$ a (partial) functional.
 Then there exists a pca $A[F]$ defined on the set $A$ such that the identity relation
 on $A$ is a decidable applicative morphism $\iota_F: A \rightarrow A[F]$ and:
 \begin{enumerate}
  \item $F \in E_{A[F]} = E_{\iota_F}$
  \item Whenever $\gamma: A \rightarrow B$ is a decidable applicative morphism such that $F \in E_\gamma$, 
  there exists a unique (up to equivalence) decidable applicative morphism $\gamma_F: A[F] \rightarrow B$ such that $\gamma_F \circ \iota_F = \gamma$.
  Moreover, if $\delta: A[F] \rightarrow B$ is a decidable applicative morphism such that $\delta \circ \iota_F \cong \gamma$, then $\delta \cong \gamma_F$.
 \end{enumerate}
 \begin{proof}
 We will construct $A[F]$ as $A[f]$ for
 \[
  f = \bigcup_{\alpha} f_\alpha,
 \]
 where $\{f_\alpha\}_\alpha$ is a compatible family of partial functions indexed by the ordinals.
 We define this family by transfinite induction.
 For all $a,b \in A$, we let:
 \begin{itemize}
  \item $f_0 = \emptyset$
  \item $f_{\alpha+1}(a) = b$ whenever there exists $g \in \text{Tot}_{A[f_\alpha]}$ such that $a \in I_1^{A[f_\alpha]}(g)$ and $F(g) = b$.
  \item For $\lambda > 0$ a limit ordinal, $f_\lambda = \bigcup_{\alpha < \lambda} f_\alpha$.
 \end{itemize}
 It is not hard to see that this indeed defines a compatible family of partial functions.
 So $f = \bigcup_{\alpha} f_\alpha$ is a partial function, and we define $A[F] := A[f]$.
 We denote the application in $A[F]$ by $\cdot^F$.

  For part (i), suppose $g \in \text{Tot}_{A[F]}$, $a \in I_1^{A[F]}(g)$, and $F(g)$ is defined.
  Then it is not hard to see that there is $\alpha$ such that $g \in \text{Tot}_{A[f_\alpha]}$, $a \in I_1^{A[f_\alpha]}(g)$.
  Then $f_{\alpha+1}(a) = F(g)$, hence $f(a) = F(g)$.
  Indeed, $F$ is an effective operation in $A[F]$.

For part (ii):
Suppose $\gamma: A \rightarrow B$ is a decidable applicative morphism such that $F \in E_\gamma$.
Let $r_F \in I_2^\gamma(F)$.
  Denote the realizer for $\gamma$ by $r$.
  Let $\mathsf{d}$ be the decider for $\gamma$.
  Let $\mathsf{c}$ be such that for all $x \in \gamma(a)$, $v \in \gamma(u)$, $\mathsf{c}xy \in \gamma([a]*u)$.
  Let $\mathsf{c}'$ be such that for all $y \in \gamma(a), x \in \gamma(u)$, $\mathsf{c}'xy \in \gamma (u*[a])$.
  Let $\mathsf{q}_0,\mathsf{q}_1$ be such that for all $x \in \gamma(a)$, $\mathsf{q}_0 x \in \gamma(\mathsf{p}_0a)$, $\mathsf{q}_1 x \in \gamma(\mathsf{p}_1a)$.
  By the recursion theorem, let $e \in B$ be such that
  \begin{align}
   exv \lesssim \textsf{ if } &\mathsf{d}\mathsf{q_0}(re(\mathsf{c}xv)) \textsf{ then } \mathsf{q_1}(re(\mathsf{c}xv)) \\
   &\textsf{ else } ex(\mathsf{c}'v(r_F(r(\mathsf{q_1}(re(\mathsf{c}xv)))))).
  \end{align}
  
  Then one can verify that for $v \in \gamma([])$, we have for all $e \in \gamma(a)$, $x \in \gamma(a')$ such that $a\cdot^f a'\halts$:
  \[
   exv \in \gamma(a \cdot^f a').
  \]
  Essential is that if $a \in I_1^{A[F]}(f)$, then for $b \in \gamma(a)$, $rb \in I_1^\gamma(f)$.
   
   So by the above, $\la x y \ra xyv$ realizes $\gamma$ as an applicative morphism $\gamma_F: A[F] \rightarrow B$,
   and the fact that $\gamma = \gamma_F \circ \iota_F$ is obvious.
   
   Uniqueness is easy to check.
 \end{proof}
\end{theorem}

\begin{remark}\label{rem:stages}
 The construction of $A[F]$ can be thought of as a construction in stages: we say that $a\cdot^F b\halts$ \emph{at stage} $\alpha$
 if $a \cdot^{f_\alpha} b \halts$.
 In proofs of statements about $A[F]$, we will often use induction on the stage at which a term denotes.
 The infinitary nature of computation relative in a functional, as remarked by Kleene in \cite{KLEENE-1959}, is reflected by this construction in (transfinitely many) stages.
\end{remark}

\begin{example} We again take a look at the functional $E$ from example \ref{ex:functionalE}.
For any pca $A$, $E$ is an effective operation in $A[E]$.
As a consequence, the computable predicates in $A[E]$ are closed under existential quantification: suppose
$p \in A[E]$ is a total function in 2 variables, taking values in $\{\true,\false\}$.
Then:
\[
 Q(n) := (\exists m) p n m
\]
is clearly also a recursive predicate in $A[E]$, so it has some index $q \in A[E]$.

For $A = \mathcal{K}_1$, this produces a familiar (see \cite{KLEENE-1959}) result.
Namely, it follows (by taking complements) that every arithmetical subset of $\N$ is computable in $\mathcal{K}_1[E]$.
We will see below that the total functions in $\mathcal{K}_1[E]$ are in fact precisely the hyperarithmetical functions.
\end{example}

\begin{remark}
 One may want to define, for functionals $F,G: A^A \rightarrow A$, a pca $A[F,G]$ in which $F,G$ are both effective operations.
 It will not work to define this as $A[F][G]$. However, if one codes any tuple of functionals into one functional $H$ (e.g. $H(g) = \p F(g) G(g)$), then 
 $A[H]$ will satisfy the required properties and the analogue of theorem \ref{thm:general_recursion} for such a tuple holds.
\end{remark}

\section{Effective Operations and Kleene computability}
In \cite{KLEENE-1959}, Kleene developed a recursion theory for functionals, which included a notion of computability relative to a functional.
A short overview of this theory, as well as comparisons to other notions of higher-type computation can be found in \cite{LONGLEY-HIGHER_TYPE_COMP}. 

In this section we will show that a total function $f: \N \rightarrow \N$ is computable relative to a functional $F: \N^\N \rightarrow \N$
(in Kleene's sense) if and only if it has an index in $\mathcal{K}_1[F]$.
This result will be almost immediate, since (as we will see) Kleene's indexed set of functions computable in $F$ is basically an explicit model
of $\mathcal{K}_1[F]$.
Our construction has the advantage of being independent of any explicit model of $\mathcal{K}_1$
and avoids a lot of technicalities. Moreover, there is the immediate generalization to arbitrary pcas and the connection to realizability.

In his paper, Kleene defined an $\N$-indexed set of functions that may take tuples of arguments of all finite types.
The definition consists of nine inductive rules, named S1-S9.
A function with index $e$ is denoted by $\{e\}$.
A function $g: \N \rightarrow \N$ is then computable \emph{relative to} the functional $F: \N^\N \rightarrow \N$ if there is an index $e \in \N$
such that for all $n \in \N$:
\[
 \{e\}(n,F) = g(n).
\]
Instead of defining Kleene's S1-S9 in full generality, we will just state eight derived rules that define the indexed set of functions computable
relative to some functional $F$ according to Kleene.
The rule S7 is missing since that applies specifically to type 1 arguments, and we only consider a specific type 2 oracle.
\begin{definition}\label{dfn:kleene_recursion}
 Let $F: \N^\N \rightarrow \N$ be a functional.
 Also fix some primitive recursive coding of tuples $\la - \ra: \N^{<\omega} \rightarrow \N$.
 We define an indexed set of partial functions $\N^{<\omega} \rightarrow \N$ by the following rules:
 for all $n,m \in \N$:
 \begin{itemize}
  \item[S1] Successor function: $\{\la 1,1 \ra\}(n) = n+1$
  \item[S2] Constant function: $\{\la 2,1, m \ra\}(n) = m$ 
  \item[S3] Projection: for all $n_1,\ldots,n_m \in \N$, $\{\la 3, m \ra \}(n_1,\ldots,n_m) = n_1$
  \item[S4] Composition: For $g, h \in \N$, all $k,n_1,\ldots,n_k \in \N$:
  \[
   \{\la 4,k, g, h \ra\}(n_1,\ldots,n_k) \simeq \{g\}(\{h\}(n_1,\ldots,n_k),n_1,\ldots,n_k).
  \] 
  \item[S5] Primitive recursion: For any $g,h \in \N$, if $m = \la 5,2,g,h \ra$ then for all $k$:
\begin{align}
 \{m\}(0,n) \simeq \{g\}(n) \\
 \{m\}(k+1,n) \simeq \{h\}(k,\{m\}(k,n),n)
\end{align}
\item[S6] Permutation of arguments: For any $g \in \N, 1 \leq k < r$ and $n_1,\ldots,n_r \in \N$:
\[
 \{\la 6,r,k,g \ra\}(n_1,\ldots,n_r) \simeq \{g\}(n_{k+1},n_1,\ldots,n_k,n_{k+2},\ldots,n_r).
\]
\item[S8] Application of $F$: For every $h,e \in \N$:
\[
 \{\la 8, 1,h\ra \} (n) \simeq F(\lambda k. \{h\}(k,n))
\]
\item[S9] Index vocation: for all $m,k,l,n_1,\ldots,n_{k+l} \in \N$
\[
 \{\la 9,k,l\ra\}(m,n_1,\ldots,n_{k+l}) \simeq \{m\}(n_1,\ldots,n_k).
\]
 \end{itemize}
 \begin{note}
	 Here, $a \simeq b$ has the same meaning as before: the right-hand side is defined if and only if the left-hand side is defined, and when this is the case both sides are equal.
 \end{note}
\end{definition}

In what follows, we assume that we have fixed some functional $F: \N^\N \rightarrow \N$,
and defined the partial map $\{-\}(-): \N \times \N \rightarrow \N$ according to the above.
\begin{theorem}
 The partial map $\{-\}(-): \N \times \N \rightarrow \N$ defines a partial combinatory algebra $\mathcal{K}_1^F$.
 \begin{proof}
  Tedious exercise using the clauses in definition \ref{dfn:kleene_recursion}.
 \end{proof}
\end{theorem}

The following can be viewed as the type 2 analogue of corollary \ref{cor:oosten_original}(iii).
\begin{theorem}\label{thm:kleene_iso}
There is an isomorphism of pcas $\mathcal{K}_1[F] \cong \mathcal{K}_1^F$.
\begin{proof}
 By Church's thesis, the identity relation on $\N$ is an applicative morphism $\gamma: \mathcal{K}_1 \rightarrow \mathcal{K}_1^F$.
 It is easy to see that $F \in E_\gamma$.
 By theorem \ref{thm:force_eff_operation}, the identity relation on $\N$ is an applicative morphism $\gamma_F: \mathcal{K}_1[F] \rightarrow \mathcal{K}_1^F$.
 
 It is not hard to show that we can define the application in $\mathcal{K}_1^F$
 computably in $\mathcal{K}_1[F]$ using the clauses \ref{dfn:kleene_recursion} and the recursion theorem when needed.
 So the identity relation on $\N$ is also an applicative morphism $\mathcal{K}_1^F \rightarrow \mathcal{K}_1[F]$.
 Therefore $\mathcal{K}_1[F] \cong \mathcal{K}_1^F$.
 \end{proof}
\end{theorem}

\begin{corollary}\label{cor:functionalE}
 Let $E$ be as in example \ref{ex:functionalE} for $A = \mathcal{K}_1$.
 The total functions in $\mathcal{K}_1[E]$ are precisely the hyperarithmetical functions.
 \begin{proof}
  This follows from the same fact for $\mathcal{K}_1^E$ and theorem \ref{thm:kleene_iso}, which is theorem XLVIII in \cite{KLEENE-1959}.
 \end{proof}
\end{corollary}
The above corollary is what we can deduce from Kleene's original paper on recursive functionals. However, we can do a little better using the results of Hinman's approach in \cite{HINMAN}.
There, the domains of the $\N-$indexed partial functions $\N\rightarrow\N$ ``computable in $E$'' are precisely
the $\Pi^1_1$-sets. Here we can prove something similar. In the context of pcas and applicative morphisms, domains of partial recursive functions are not particularly well-behaved (we have to impose more conditions on morphisms to preserve this structure, see for example the discussion of dominances in \cite{LONGLEY}).
Theorem \ref{thm:pca_iso_strict_pca} also exemplifies this fact.
We therefore avoid domains, and show the following instead:
\begin{proposition}\label{prop:pi_1_1_complete}
	The application in $\mathcal{K}_1[E]$ is $\Pi^1_1$-complete.
 \begin{proof}
 We show here that the application function is $\Pi^1_1$, since that does not follow immediately from Hinman's results.
 Let $r \in \N$ be the uniform index for $f$ in $\mathcal{K}_1[f]$, for any $f$.
 Recall that $\mathcal{K}_1[E] = \mathcal{K}_1[g]$, where we defined $g$ in stages.
 I'll show that $g$ is $\Pi^1_1$.
 Let $e$ be an index such that whenever $\phi$ is the characteristic function of the graph of a partial function $f: A \rightarrow A$,
 then for all $b,x \in \mathcal{K}_1$:
 \[
  eb \cdot^\phi x \simeq b \cdot^f x.
 \]
Convince yourself that such $e$ exists; since $\N$ is enumerable, we can just ``wait'' for values of $f$ to appear in the graph.
 
 Define the following predicate on $\N \times \N^\N$:
 \begin{align}\label{eq:uberpredicate}
  \Phi([b,i],X) = &(\forall s \exists t \forall y):
  \left( eb \cdot^X s \halts \land 
  (eb \cdot^X t = 0 \lor i = 1) \land (eb\cdot^X y \neq 0 \lor i = 0).\right)
 \end{align}
 It is not hard to see that $\Phi$ is $\Pi^1_1$ (in fact it is arithmetical) (to work out the details, one needs to use
 Kleene's $T-$predicate and existential quantification over dialogues).
 Let $G$ be the characteristic function of the graph of $g$.
 Now I claim:
 \begin{equation}\label{eq:Pi11_predicate}
  G([b,i]) \iff (\forall X)(\exists p,j)(\Phi([p,j],X) \Rightarrow X([p,j])) \Rightarrow X([b,i]).
 \end{equation}
\begin{proof}[Proof of claim]
First, observe:
\begin{align}
 \Phi([b,i],G) &\Leftrightarrow (\forall s)(\exists t \forall y) eb \cdot^G s \halts \land (eb \cdot^G t = 0 \lor i = 1) \land (eb\cdot^G y \neq 0 \lor i = 0) \\
 &\Leftrightarrow (\forall s)(\exists t \forall y) b \cdot^g s \halts \land (b \cdot^g t = 0 \lor i = 1) \land (b\cdot^g y \neq 0 \lor i = 0)  \\
 &\Leftrightarrow G([b,i]).
\end{align}
This proves the ``$\Leftarrow$'' direction of \eqref{eq:Pi11_predicate} (we substitute $G$ for $X$).

Now assume $G([b,i])$, that is:
\[
 (\forall s)(\exists t \forall y) b \cdot^g s \halts \land (b \cdot^g t = 0 \lor i = 1) \land (b\cdot^g y \neq 0 \lor i = 0).
\]
Suppose we have $X$ such that
\[
 (\forall p,j): \Phi([p,j],X) \Rightarrow X([p,j]).
\]
We need to prove that $X([b,i])$.
I'll prove this by showing that
\begin{equation}\label{eq:XextendsG}
 (\forall p,j): G([p,j]) \Rightarrow X([p,j]).
\end{equation}
This is done by induction on stages (see remark \ref{rem:stages}). Recall that $g$ is constructed as $g = \bigcup_\alpha g_\alpha$.
Denote the corresponding graphs by $G_\alpha$.
Then for $\alpha = 0$,
\[
 (\forall p,j): G([p,j]) \Rightarrow X([p,j])
\]
clearly holds, since $g_0 = \emptyset$.
If it holds up to $\alpha$, then:
\begin{align}
 G_{\alpha+1}[p,j] &\Rightarrow 
 (\forall s)(\exists t \forall y) p \cdot^{g_\alpha} s \halts \land (p \cdot^{g_\alpha} t = 0 \lor j = 1) \land (p\cdot^{g_\alpha} y \neq 0 \lor j = 0) \\
 &\Rightarrow 
 (\forall s)(\exists t \forall y) e p \cdot^{X} s \halts \land (e p \cdot^{X} t = 0 \lor j = 1) \land (ep\cdot^X y \neq 0 \lor j = 0) \\
 &\Rightarrow \Phi([p,j],X) \\
 &\Rightarrow X([p,j]).
\end{align}
 For limit ordinals $\lambda$, $G_{\lambda}[p,j] \iff G_\alpha[p,j]$ for some $\alpha < \lambda$ so that follows trivially.
 
 Hence we have shown \eqref{eq:XextendsG}, from which $X[b,i]$ follows by assumption.
\end{proof}

Since the right-hand side of \eqref{eq:Pi11_predicate} is $\Pi^1_1$, we are done.

To see that the application is $\Pi^1_1$ complete, see \cite{HINMAN}, section VI.1.
There one can find a proof in which a well-known $\Pi^1_1$-complete set is reduced to a
\emph{semi-recursive} set computable from $E$ in Hinman's sense (i.e. a domain of a partial recursive function). This proof is easily
adopted to a proof of $\Pi^1_1$-completeness for our case.
 \end{proof}

\end{proposition}
\section{Local Operators}
In \cite{OOSTEN:arXiv1301.0735}, a notion of realizability is studied that uses a specific local operator $\J$ on the Effective Topos.
We will show here that some of the results in that paper can be derived easily using the functional $E$ from example \ref{ex:functionalE} (for $A = \mathcal{K}_1$)
and theorem \ref{thm:force_eff_operation}.

We first give a definition of local operators in realizability toposes on a pca $A$, that we assume fixed from now on.
The topos-theoretic origin is omitted here, for more about that we refer to \cite{OOSTEN} or \cite{PITTS}.
We introduce the following notation: For $P,Q \subseteq A$ subsets, let
\[
 P \rightarrow Q = \{a \in A \varco (\forall b \in P) ab \halts \text{ and } ab \in Q\}
\]

\begin{definition}\label{dfn:local_operator}
 A \emph{local operator} is a function $\J: \mathcal{P}(A) \rightarrow \mathcal{P}(A)$ such that there are $e_1,e_2,e_3 \in \N$ that satisfy:
 \begin{enumerate}
  \item $e_1 \in \bigcap_{P \subseteq A} P \rightarrow \J P$
  \item $e_2 \in \bigcap_{P \subseteq A} \J \J P \rightarrow \J P$
  \item  $e_3\in \bigcap_{P, Q\subseteq A} \left(P \rightarrow Q\right) \rightarrow \left(\J P \rightarrow \J Q\right)$
 \end{enumerate}
\end{definition}

We define a pre-order on local operators as follows: we say $\J \leq \J'$ whenever there is $e \in A$ such that:
\[
 e \in \bigcap_{P \subseteq A} \J P \rightarrow \J' P.
\]
We write $\J \cong \J'$ if $\J \leq \J'$ and $\J' \leq \J$.
There is a least local operator $\J_\bot$, given by the identity: $\J_\bot(P) = P$.
There is also a biggest local operator given by $\J_\top(P) = A$.

The lattice of local operators (up to isomorphism) corresponds to the lattice of subtoposes of the realizability topos on $A$, denoted $\rt{A}$.
Some of these subtoposes are realizability toposes $\rt{B}$ themselves, in that case we have found a \emph{geometric inclusion} $\rt{B} \rightarrow \rt{A}$.
\begin{definition}[\cite{OOSTEN-HOFSTRA}]\label{dfn:computationally_dense}
 An applicative morphism $\gamma: A \rightarrow B$ is called \emph{computationally dense} if there exists $m \in B$ such that
 \[
  (\forall b \in B)(\exists a \in A)(\forall a') b\gamma(a')\halts \Rightarrow m\gamma(aa') \subseteq b\gamma(a').
 \]
\end{definition}

Geometric inclusions $\rt{B} \rightarrow \rt{A}$ correspond precisely to computationally dense applicative morphisms $\gamma: A \rightarrow B$ that satisfy the additional condition
\[
 \text{(in) } (\exists e) (\forall b) (\exists a) (eb \in \gamma(a) \text{ and } m\gamma(a) = \{b\})
\]
where $m$ is as in definition \ref{dfn:computationally_dense} (this is
a consequence of results in \cite{OOSTEN-HOFSTRA,JOHNSTONE-2013}). In that case the corresponding local operator in $\rt{A}$ is given by
\[
 \J(P) = \{a \in A \varco m\gamma(a) \subseteq \bigcup_{a' \in P} \gamma(a')\}.
\]
The following corresponds to proposition 2.2 in \cite{OOSTEN-GENERAL_RECURSION}.
\begin{proposition}
 For $f: A \rightarrow A$ a partial function, the identity relation on $A$ is a computationally dense applicative morphism $\iota_f: A \rightarrow A[f]$
 that satisfies the condition (in). Therefore, there is a canonical geometric inclusion $\rt{A[f]} \rightarrow \rt{A}$.
\end{proposition}

By our construction, we obtain for every (partial) functional $F: A^A \rightarrow A$ a subtopos $\rt{A[F]} \rightarrow \rt{A}$
corresponding to the local operator
\[
 \J_F(P) = \{a \in A \varco m\cdot^F a \in P\}.
\]
One can check that this indeed defines a local operator according to definition \ref{dfn:local_operator}.

\subsection{\texorpdfstring{$\J$-Assemblies}{J-Assemblies}}
In \cite{HYLAND}, Martin Hyland showed that for every function $f: \N \rightarrow \N$,
there is a least local operator that forces $f$ to be computable in the Effective Topos.
This yields an embedding of the Turing degrees into the lattice of local operators.
In \cite{PHOA}, Wesley Phoa showed that the corresponding subtopos is equivalent to $\rt{\mathcal{K}_1^f}$, the realizability topos on the pca of recursive application with oracle $f$.
This result can be generalized to arbitrary pcas and partial functions: for (partial) functions $f: A \rightarrow A$,  the least subtopos of $\rt{A}$ in which $f$ is realizable is equivalent to $\rt{A[f]}$. For a proof, see \cite{GEOMMORT}.

We will now carry out a similar statement for effective operations.
To make things easier, we will work on a category of assemblies.
The consequences for realizability toposes are immediate for anyone familiar with them.
In the following we assume that we have fixed some pca $A$.
\begin{definition}
 Let $\J$ be a local operator for $A$. A $\J$-\emph{assembly} is a pair $(X,E)$ where $E: X \rightarrow A$ is a total relation.
 
 For $(X,E), (Y,F)$ $\J$-assemblies, a \emph{morphism} of $\J-$assemblies is a function $f: X \rightarrow Y$ together with a $t \in A$ such that
 \[
  t \in \bigcap_{x \in X} E(x) \rightarrow \J F(f(x)).
 \]
 We say that $t$ \emph{tracks} $f$.
\end{definition}

One can check that the composition $g \circ f$ of morphisms of $\J-$assemblies $f: (X,E) \rightarrow (Y,F)$, $g: (Y,F) \rightarrow (Z,G)$ is again a morphism of $\J$-assemblies,
so that $\J-$assemblies form a category $\ass{A}_\J$.
The category $\ass{A}_\J$ has a lot of structure: it is regular, cartesian closed and has finite colimits.
For $(X,E)$, $(Y,F)$ $\J$-assemblies, the exponential $(Z,G) := (Y,F)^{X,E}$ is given by:
\begin{align}
 &Z = \{f: X \rightarrow Y \varco f \text{ is a morphism } (X,E) \rightarrow (Y,F) \} \\
 &G(f) = \{t \in A \varco t \text{ tracks } f\}.
\end{align}

The assembly $\mathbf{A} = (A, \text{id}_A)$ can be seen as a representation of the pca $A$ in the category of $\J$-assemblies.
The object $\mathbf{A}^\mathbf{A}$ consists of the morphisms $\mathbf{A} \rightarrow \mathbf{A}$ in $\ass{A}_\J$.
For $\J = \J_\bot$, these are precisely the computable functions in $A$.
Also, one can check that the arrows
\[
 \mathbf{A}^\mathbf{A} \rightarrow \mathbf{A}
\]
in $\ass{A}_{\J_\bot}$ are precisely the effective operations in $A$.
We have the following theorem:
\begin{theorem}\label{thm:llo_functional}
 Let $F: A^A \rightarrow A$ be a total functional, and assume that $F$ defines a morphism
 \[
  \mathbf{A}^\mathbf{A} \rightarrow \mathbf{A}
 \]
in $\ass{A}_\J$. Then $\J_F \leq \J$.
 \begin{proof}
 We can assume that $\J$ preserves inclusions, otherwise we let
 \[
  \J'P = \bigcup_{Q \subseteq P} \J P
 \]
and show that $\J' \cong \J$, which is an easy exercise using definition \ref{dfn:local_operator}.
 
 By computational density, there exists $t \in A$ such that for all $v,x \in A$:
 \[
  m\cdot^F (tvx) \lesssim v\cdot^F x.
 \]
 Since $m$ is independent of $F$, we also know that if the right side halts at some stage, the left hand side also halts at that stage.
 
We pick $e_1,e_2,e_3$ as in definition \ref{dfn:local_operator}. 

Suppose that $F: \mathbf{A}^\mathbf{A} \rightarrow \mathbf{A}$ is tracked by $r$ as morphism of $\J$-assemblies.
 
 Choose $\mathsf{m} \in A$ (using the recursion theorem) so that:
 \begin{align}
	 \mathsf{m}xy &\lesssim \textsf{if } \p_0 (m([x]*y)) \\
	 &\textsf{then } e_1(\p_1(m([x]*y))) \\
	 &\textsf{else } e_2 \left(e_3(\la z \ra \mathsf{m}x(y*[z])) r(\la w \ra (\mathsf{m}(tvw)[]))\right) \\
	 &\text{ where } v = p_1(m([x]*y)).
 \end{align}
 
 Here $[]$ is the empty sequence.
 We will prove by induction on the stage at which $m\cdot^F x\halts$ that for all $x$,
 \begin{equation}\label{eq:inductiongoal}
	 \mathsf{m}x[] \in \J\{m \cdot^F x\}.
 \end{equation}
 
 Stage $0$: in that case $m([x]) = \p \true c$ for some $c$.
 Then
\begin{align}
	\mathsf{m}x[] = e_1(\p_1(m([x]*[]))) \in \J\{c\}.
\end{align}

Stage $\lambda > 0$: Let $u = [u_0,\ldots,u_n]$ be the halting dialogue between $m$ and $x$.
It is easy to see that
\[
	\mathsf{m}x[u_0,\dots,u_n] \in \J\{m\cdot^F x\}.
\]

Now suppose that $0 < i \leq n+1$ is such that
\begin{equation}\label{eq:assumption_reverse_ind}
	\mathsf{m}([x]*u^{<i+1}] \in \J\{m \cdot^F x\}.
\end{equation}
We know that:
\[
 m([x]*u^{<i}) = \p \false v
\]
where $v$ is such that at some stage $\alpha < \lambda$,
for each $w$, $v \cdot^F w \halts$.
So for each $w$, also $m\cdot^F (tvw) \halts$ at stage $\alpha$.
Let $f$ be the function $w \mapsto v \cdot^F w$.
By induction hypothesis,
\[
	\mathsf{m}(tvw)[] \in  \J\{v \cdot^F w\}
\]
for each $w$, so $\la w \ra \mathsf{m}(tv w)[]$ tracks $f$ as morphism of $\J-$assemblies.
So \[r(\la w \ra (\mathsf{m}(tv w)l)) \in \J\{F(f)\}\]
and $F(f) = u_i$.
By \eqref{eq:assumption_reverse_ind}, we have:
\[
	\la z \ra \mathsf{m}([x]*u^{<i}*[z]) \in \{F(f)\} \rightarrow \J\{m \cdot^F x\}
\]
and therefore
\[
	e_2 \left(e_3(\la z \ra \mathsf{m}x(u^{<i}*[z])) r(\la w \ra (\mathsf{m}(tvw)[]))\right) \in \J\{m \cdot^F x\}
\]
hence it follows that
\begin{equation}
	\mathsf{m}([x]*u^{<i}) \in \J\{m \cdot^F x\}.
\end{equation}

By reverse induction, it follows that
\[
	\mathsf{m}x [] \in \J\{m \cdot^F x\},
\]
so this finishes the induction step, so we have \eqref{eq:inductiongoal}.

Therefore
\[
	\la x \ra \mathsf{m}x [] \in \bigcap_{P \subseteq A} \J_F P \rightarrow \J P
\]
and we are done.
 \end{proof}
\end{theorem}

\subsection{Pitts' local operator}
In this section we will take a look at the operator $\J$ that is studied in \cite{OOSTEN:arXiv1301.0735}.
This operator was introduced by Andrew Pitts in his thesis (example 5.8, \cite{PITTS}).
Denote by $\nabla(\N)$ the assembly
\[
 (\N, T) \text{ where } T(n) = \N \text{ for all }n.
\]
We can define the following subobject of $\nabla(\N)$:
\[
 (\N, R) \rightarrowtail \nabla(\N) \text{ where } R(n) = \{m \in \N \varco m \geq n \}.
\]
Pitts' local operator $\J$ is defined as the least local operator in $\rt{\mathcal{K}_1}$ that forces the above arrow to be an isomorphism of $\J-$assemblies.
Explicitly, it is given by:
\[
 \J(P) = \bigcap \{ Q \subseteq \N \varco \{0\} \land P \subseteq Q \text{ and } \{1\} \land (\bigcup_{n\in \N} (R(n) \rightarrow Q)) \subseteq Q \}
\]
where for $R,S \subseteq \N$:
\[
 R \land S = \{\p rs \varco r \in R, s \in S\}.
\]
Write $\mathbf{N}$ for the assembly $\mathbf{A}$ for $A = \mathcal{K}_1$, so
\[
 (\N,N) \text{ where } N(n) = \{n\}.
\]
Recall the functional $E$ from example \ref{ex:functionalE}.
For $A = \mathcal{K}_1$, we rather define it as
\[
 E(f) = \begin{cases}
         0 &\text{ if } (\exists n) f(n) \\
         1 &\text{ otherwise}.
        \end{cases}
\]

We have the following proposition:
\begin{proposition}\label{prop:J1hyperarith}
 There is an index $r_E \in \N$ that tracks $E$ as morphism of $\J$-assemblies ${\mathbf{N}}^{\mathbf{N}} \rightarrow \mathbf{N}$.
\begin{proof}
	Pick $e_3$ as in definition \ref{dfn:local_operator}. 

 By corollary 1.6 of \cite{OOSTEN:arXiv1301.0735}, there is a partial recursive function $G$ such that for
 all \[x_0 \in \J\{a_0\}, \ldots, x_{n-1} \in \J\{a_{n-1}\},\] we have
 \begin{align}
  G(\la x_0, \ldots,x_{n-1} \ra) \in \J\{0\} &\text{ if for some }i<n, a_i = 0 \\
  G(\la x_0, \ldots, x_{n-1} \ra) \in \J\{1\} &\text{ otherwise }
 \end{align}
Now let $t_E$ be an index in $\mathcal{K}_1$ defined by:
\[
 t_E e n \simeq G(\la e (0), \ldots, e (n)\ra)
\]
Then whenever $e$ realizes a total function $f: \mathbf{N} \rightarrow \mathbf{N}$, we $(\exists n) f(n) = 0$ if and only if $(\exists n) e n \in \J\{0\}$,
which holds if and only if there is $n$ such that
\[
 t_E e \in \{m \varco m \geq n\} \rightarrow \J\{0\}
\]
since otherwise
\[
 t_E e \in \{m \varco m \geq 0\} \rightarrow \J\{1\}.
\]

By definition of $\J$, we have in the first case that $t_E e\in \J\{0\}$, and in the latter case $t_E e \in \J\{1\}$.
Now $t_E$ tracks $E$ as morphism of $\J-$assemblies ${\mathbf{N}}^{\mathbf{N}} \rightarrow \mathbf{N}$.
\end{proof}
\end{proposition}

The following corollary corresponds to theorem 2.2 in \cite{OOSTEN:arXiv1301.0735}.
\begin{corollary}\label{cor:J_hyperarithmetical}
 The morphisms of $\J$-assemblies $\mathbf{N} \rightarrow \mathbf{N}$ are precisely the hyperarithmetical functions.
 \begin{proof}
 By theorem \ref{thm:llo_functional}, $\J_E \leq \J$.
 Therefore it is easily seen that every morphism of $\J_E$ assemblies $\mathbf{N} \rightarrow \mathbf{N}$
 is a morphism of $\J-$assemblies. Now it follows that every hyperarithmetical function is a morphism of $\J$-assemblies by corollary \ref{cor:functionalE}.
 
  For the converse, observe that every morphism $f: \mathbf{N} \rightarrow \mathbf{N}$ is $\Pi^1_1$ since:
  \[
   f(n) = m \iff (\forall B) \{\p 0 m \} \subseteq B \land \{1\}\land (\{e \varco (\exists n)(\forall m \geq n) e m \in B\}) \subseteq B \rightarrow e(\mathsf{a}n) \in B
  \]
where $e$ is a realizer for $f$ as a morphism $\mathbf{N} \rightarrow \mathbf{N}$.
Notice that the above expression is $\Pi^1_1$.
Since $f$ is total, it is also $\Sigma^1_1$.

Therefore the morphisms of $\J$-assemblies $\mathbf{N} \rightarrow \mathbf{N}$ are precisely the hyperarithmetical functions.
 \end{proof}
\end{corollary}

\subsection{\texorpdfstring{$\J$-representable functions}{J-representable functions}}
In \cite{OOSTEN:arXiv1301.0735}, the morphisms of $\J$-assemblies $\mathbf{N} \rightarrow \mathbf{N}$
are called $\J-$\emph{representable}. Since for every $n, m \in \N$
\begin{equation}
 n \neq m \iff \J\{n\} \cap \J\{m\} = \emptyset,
\end{equation}
the following is well defined:
\begin{definition}
 A partial function $F: \N \rightarrow \N$ is $\J-$\emph{representable} if there exists $e \in \N$ such that for all $m \in \N$:
 \[
  \phi_e(n) \in \J\{m\} \iff F(n) = m.
 \]
Here $\phi_e$ denotes the $e-$th partial computable function from ordinary recursion theory (or the application in $\mathcal{K}_1$).

By defining $\psi_e = F$ for such $e$, we obtain an indexed set of partial functions and a new application
\[
 e \cdot n \mapsto \psi_e(n)
\]
on $\N$.
\end{definition}

A natural question is whether the application $e \cdot n \simeq \psi_e(n)$ is combinatory complete in the sense of definition \ref{dfn:pca}. This turns out to be the case:
\begin{proposition}\label{prop:J_1_pca}
 The application $e\cdot n \simeq \psi_e(n)$ defines a pca $\J_1$ on $\N$.
 \begin{proof}
  This is an easy exercise on pcas and definition \ref{dfn:local_operator}.
 \end{proof}
\end{proposition}

The above proposition would be particularly hard to prove (if not false) if one uses the strict version of the definition of a pca. In the next section we will see that we do not have to worry about this.

\subsection{Strict Pcas} As announced in the introduction, several authors use a stronger notion of combinatory completess, that we shall call \emph{strictly combinatory complete}:
\begin{definition}\label{dfn:strict_pca}
 Let $A$ be a partial applicative structure. Then $A$ is called a \emph{strict pca} when there are $\k,\s \in A$ such that for all $a,b,c \in A$:
 \begin{enumerate}
  \item $\s a b \halts$ 
  \item $\k a b = a$
  \item $\s a b c \simeq a c (b c)$.
 \end{enumerate}
We also call the applicative structure on $A$ \emph{strictly combinatory complete}.
\end{definition}
This definition might be preferable in some cases, for example when one would like that the set of domains
of partial computable functions in a pca $A$ form a \emph{dominance} in $\rt{A}$(see \cite{LONGLEY}).
However, as mentioned before, these domains are not part of the structure of a pca in the category PCA.

The following theorem supports this: in the context of applicative morphisms, the notions of a pca and a strict pca are essentially the same.
\begin{theorem}\label{thm:pca_iso_strict_pca}
 Every pca is isomorphic to a strict pca.
 \begin{proof}
  Let $A$ be a pca. If $A$ is trivial, then the statement is obvious, so we assume $A$ non-trivial.
  
  The idea is to define an explicit pca structure $A'$ on the set $A$, such $a \mapsto \{a\}$ is an applicative morphism
  $A' \rightarrow A$. In this explicit structure, we also make sure that we can compute the application in $A$.

  To define the new applicative structure on $A$, we use tuples and Curry numerals from $A$ to ensure that what we define is actually computable in $A$.
  We denote the new application by $\cdot: A \times A \rightarrow A$, and the resulting partial applicative structure by $A'$.
  The definition is by induction on the following clauses:
  \begin{enumerate}
  \item Constant function: for all $a,b, \in A$:
  \begin{align}
   [\overline{1}] \cdot a = [\overline{1},a] \\
   [\overline{1},a] \cdot b = a
  \end{align}
\item S-combinator: For $a,b,c \in A$:
\[
 [\overline{2}]\cdot [a,b,c] \simeq a\cdot c \cdot (b\cdot c)
\]
\item Currying of parameters: for all $e,a,b,c \in A$:
\begin{align}
 &[\overline{3},e] \cdot a = [\overline{3},e,a] \\
 &[\overline{3},e,a] \cdot b = [\overline{3},e,a,b] \\
 &[\overline{3},e,a,b]\cdot c \simeq e \cdot [a,b,c] 
\end{align}
\item Application in $A$:
	\begin{align}
		&[\overline{4}] \cdot a = [\overline{4},a] \\
		&[\overline{4},a] \cdot b \simeq ab
	\end{align}
\end{enumerate}

Since $A$ is non-trivial, the above is well-defined.
Strict combinatory completeness of $A$ is satisfied by taking $\k' = [ \overline{1} ]$, $\s' = [\overline{3},[\overline{2}]]$. It is easy to see that $a \mapsto \{a\}$ can be realized as applicative morphism $A' \rightarrow A$,
and as applicative morphism $A \rightarrow A'$ it is realized by $[\overline{4}]$.
\end{proof}
\end{theorem}
\begin{remark}
Note that the proof above is entirely constructive; to define the application it only uses natural induction
and explicit constructions with the combinators $\k$ and $\s$.

Another observation is that if we would define $A'$ using only the clauses 1-3, we obtain a strict pca
$A'$, together with an applicative morphism $A' \rightarrow A$ given by $a \mapsto \{a\}$, that has ``forgotten'' a lot of the structure of $A$. We can then adjoin certain partial functions in $A$ to $A'$ to get
some of them back (such as the application, so that we obtain $A$ again).
In the case of a decidable pca, for instance, we would lose decidability in this way. So, not only can we adjoin functions to a pca as an oracle, we are also able to forget (all but some) functions in a pca.
\end{remark}

Finally, we have the following theorem:
\begin{theorem}  The pca $\J_1$ is isomorphic to $\mathcal{K}_1[E]$.
 \begin{proof}
  Let $\gamma: \mathcal{K}_1 \rightarrow \J_1$ be the applicative morphism $n \mapsto \{n\}$ (it is easily verified that this is a decidable applicative morphism).
  
  Now $E \in E_\gamma$ by proposition \ref{prop:J1hyperarith}.
  Therefore, by theorem \ref{thm:force_eff_operation}, $\gamma$ factors through $\mathcal{K}_1[E]$, hence we have a decidable applicative morphism
  $\gamma_E: \mathcal{K}_1[E] \rightarrow \J_1$, given by the identiy on $\N$.
  
  This morphism also goes the other way: this follows from the fact that the application in $\mathcal{K}_1[E]$ is $\Pi^1_1$-complete (proposition \ref{prop:pi_1_1_complete}), and the fact that
  the application in $\J_1$ is $\Pi^1_1$.
 \end{proof}
\end{theorem}

\section{Higher types}
One may wonder whether theorem \ref{thm:kleene_iso} also holds for computability in higher type objects (so, type 3 and higher).
The answer to that question is positive: for functionals $x_1,\ldots,x_n$ of certain types, one can define in a similar fashion a pca
\[
 \mathcal{K}_1[x_1,\ldots,x_n]
\]
that is again isomorphic to the pca of recursive application relative relative to $x_1,\ldots,x_n$ in Kleene's sense.
The idea is to define $\mathcal{K}_1[x_1,\ldots,x_n] = \mathcal{K}_1[f]$, where $f$ is defined in stages, simultaneously for 
$x_1,\ldots,x_n$, so that the resulting applicative structure is closed under application of any of the functionals.

However, the connection with effective operations does not generalize in this way.
Furthermore, it is not known (at least to me) what an analogue of theorem \ref{thm:llo_functional} would be for, say, a type 3 functional. This might be interesting to find out.

\bibliographystyle{plain}

\begin{thebibliography}{10}
\bibitem{GEOMMORT}
Eric Faber and Jaap van Oosten.
\newblock {M}ore on {G}eometric {M}orphisms between {R}ealizability {T}oposes.
\newblock 2014.
\newblock Submitted.

\bibitem{FREY}
Jonas Frey.
\newblock A characterization of realizability toposes.
\newblock 2014.
\newblock Submitted to {J}ournal of {P}ure and {A}pplied {A}lgebra.

\bibitem{HYLAND-GANDY}
R.~O. Gandy and J.~M.~E. Hyland.
\newblock Computable and recursively countable functions of higher type.
\newblock In {\em Logic {C}olloquium 76 ({O}xford, 1976)}, pages 407--438.
  Studies in Logic and Found. Math., Vol. 87. North-Holland, Amsterdam, 1977.

\bibitem{HINMAN}
Peter~G. Hinman.
\newblock {\em Recursion-theoretic hierarchies}.
\newblock Springer-Verlag, Berlin, 1978.
\newblock Perspectives in Mathematical Logic.

\bibitem{OOSTEN-HOFSTRA}
Pieter Hofstra and Jaap van Oosten.
\newblock Ordered partial combinatory algebras.
\newblock {\em Math. Proc. Cambridge Philos. Soc.}, 134(3):445--463, 2003.

\bibitem{HYLAND}
J.~M.~E. Hyland.
\newblock The effective topos.
\newblock In {\em The {L}.{E}.{J}. {B}rouwer {C}entenary {S}ymposium
  ({N}oordwijkerhout, 1981)}, volume 110 of {\em Stud. Logic Foundations
  Math.}, pages 165--216. North-Holland, Amsterdam, 1982.

\bibitem{JOHNSTONE-2013}
Peter Johnstone.
\newblock Geometric morphisms of realizability toposes.
\newblock {\em Theory Appl. Categ.}, 28:No. 9, 241--249, 2013.

\bibitem{KLEENE-1959}
S.~C. Kleene.
\newblock Recursive functionals and quantifiers of finite types. {I}.
\newblock {\em Trans. Amer. Math. Soc.}, 91:1--52, 1959.

\bibitem{LONGLEY}
John~R. Longley.
\newblock {\em Realizability Toposes and Language Semantics}.
\newblock PhD thesis, University of Edinburgh, 1995.

\bibitem{LONGLEY-HIGHER_TYPE_COMP}
John~R. Longley.
\newblock Notions of computability at higher types. {I}.
\newblock In {\em Logic {C}olloquium 2000}, volume~19 of {\em Lect. Notes
  Log.}, pages 32--142. Assoc. Symbol. Logic, Urbana, IL, 2005.

\bibitem{PHOA}
Wesley Phoa.
\newblock Relative computability in the effective topos.
\newblock {\em Math. Proc. Cambridge Philos. Soc.}, 106(3):419--422, 1989.

\bibitem{PITTS}
A.M. Pitts.
\newblock {\em The Theory of Triposes}.
\newblock PhD thesis, Cambridge University, 1981.

\bibitem{OOSTEN-GENERAL_RECURSION}
Jaap van Oosten.
\newblock A general form of relative recursion.
\newblock {\em Notre Dame J. Formal Logic}, 47(3):311--318 (electronic), 2006.

\bibitem{OOSTEN}
Jaap van Oosten.
\newblock {\em Realizability: an introduction to its categorical side}, volume
  152 of {\em Studies in Logic and the Foundations of Mathematics}.
\newblock Elsevier B. V., Amsterdam, 2008.

\bibitem{OOSTEN:arXiv1301.0735}
Jaap van Oosten.
\newblock Realizability with a local operator of {A}.{M}. {P}itts.
\newblock {\em Theoret. Comput. Sci.}, 546:237--243, 2014.

\end{thebibliography}

\end{document}